# Existence of the signal in the signal plus background model

## Tonglin Zhang[1]

*Purdue University*

**Abstract:** Searching for evidence of neutrino oscillations is an important problem in particle physics. Suppose that evidence for neutrino oscillations from an LSND experiment reports a significant positive oscillation probability, but that the LSND result is not confirmed by other experiments. In statistics, such a problem can be proposed as the detection of signal events in the Poisson signal plus background model. Suppose that an observed count $X$ is of the form $X = B + S$, where the background $B$ and the signal $S$ are independent Poisson random variables with parameters $b$ and $\theta$ respectively, $b$ is known but $\theta$ is not. Some recent articles have suggested conditioning on the observed bound for $B$; that is, if $X = n$ is observed, the suggestion is to base the inference on the conditional distribution of $X$ given $B \leq n$. This suggestion is used here to derive an estimator of the probability of the existence of the signal event. The estimator is examined from the view of decision theory and is shown to be admissible.

## 1. Introduction

In some problems, a signal $S$ may be combined with a background $B$, leaving an observed count $X = B + S$. Here we suppose that $B$ and $S$ are independent Poisson random variables with means $b$ and $\theta$ respectively so that $X$ has a Poisson distribution with mean $b + \theta$. Further $b$ is assumed known but $\theta$ is not, as might be appropriate if there were historical data on the background only. Models of this nature arise in astronomy and high energy physics. One of the most interesting problems is the issue in experiments of neutrino oscillations. Three types of neutrinos have been identified. They are denoted by $\nu_e$, $\nu_\mu$ and $\nu_\tau$ respectively. A task of the experiments is a search for $\nu_\mu$ oscillating into $\nu_e$. There are 3 sources of background. One is from the real $\nu_e$ and the other two are from two different mistaken identifications [1, 12]. For example, the KARMEN Group has been searching for neutrino oscillations reported from an earlier experiment at Rutherford Laboratory Detector in the United Kingdom. They had expected to see about 15.8 background events and had observed 15 events total [2]. This example and others have sparked interest in statistical inference when maximum likelihood estimators are on or near a physical boundary of the parameter space. Recent work along these lines is reviewed by Mandelkern [9] and discussants. Here, we consider the signal existence problem in an experiment. It could be either an estimation problem of the probability of the existence of the signal events or a testing problem of the existence of the signal events. The formulation for the estimation problem, estimating $P_\theta[S > 0|X = n]$ instead of $\theta$ itself, is equivalent to estimating the parameter $p = p_\theta(n) = b^n/(b+\theta)^n$ given observed $X = n$ since $P_\theta[S = 0|X = n] = b^n/(b+\theta)^n$.

[1]Department of Statistics, Purdue University, 150 North University Street, West Lafayette, IN 47907-2067, e-mail: tlzhang@stat.purdue.edu

The quantity $1 - p$ is the conditional oscillation probability in experiments in searching the evidence of neutrino oscillations in particle physics. In the past few years, tremendous progress has been achieved to firmly establish the nature of neutrino oscillations using neutrinos from the sun [6]. However, this situation is still unsettled in accelerator based experiments. The only evidence for oscillations in the appearance mode is from the LSND experiment [1], which reports an oscillation probability of 0.264% plus or minus measurement errors. However, the result from the LSND experiment is not confirmed from the KARMEN group [2, 7], which reports no evidence of neutrino oscillations. The result from the LSND experiment is also partially excluded from accelerator based experiments, such as the NuTev group [3]. Recently, the Mini-Boone experiment has started receiving data at Fermilab and will either confirm or refute the LSND effect, but the results won't be available until in a few more years [12]. If the neutrino oscillation results were confirmed, experiments would provide a precision measurement of the oscillation parameters [1].

Suppose that $X = n$ is observed and consider the problem of estimating $p = b^n/(b+\theta)^n$ for $\theta \geq 0$. Let $x \vee y = \max(x,y)$ and $x \wedge y = \min(x,y)$ for any $x$ and $y$. Substituting the maximum likelihood estimator (MLE) $\hat{\theta} = \max(0, n-b) = 0 \vee (n-b)$ into the expression, the MLE of $p$ is given by $p_{\hat{\theta}} = p_{\hat{\theta}}(n) = b^n/(b \vee n)^n$. That is the MLE takes the ratio of the $n$-th power of the background parameter and the $n$-th power of the maximum of the background parameter and the observed counts. An obvious disadvantage of $p_{\hat{\theta}}$ is that it is very spiky at $b$. For example, $p_{\hat{\theta}} = 1$ if $n \leq b$ and $p_{\hat{\theta}} \approx e^{-(n-b)}$ if $n$ is slightly larger than $b$ for large $n$. This indicates that only a few more observations will make the MLE significantly change from 1 to a number close to 0 no matter how large $b$ is. Note that the background $B$ is ancillary since its distribution is completely known. Even though it is not observed, its bound is. If $X = n$, then $B \leq n$. In this spirit, let $\hat{p}(n) = P(B = n|B \leq n)$, and observe that $\hat{p}$ is computable since $b$ is known. Thus, $\hat{p}$ is the conditional probability of the number of background events equal to the observed counts, given the observed bound $n \leq b$; and $0 < \hat{p} < 1$ for all $n \geq 1$. For the data reported by the KARMEN group, $b = 15.8$, $X = 15$. Since $P[B = 15|B \leq 15] = 0.206$ and $\hat{\theta} = 0$, $\hat{p} = 0.206$ and $p_{\hat{\theta}} = 1$.

The admissibility or the inadmissibility of the two estimators $\hat{p}$ and $p_{\hat{\theta}}$ will be discussed in Section 2. It is shown that $\hat{p}$ is an admissible estimator of $p$ under the squared error loss and $p_{\hat{\theta}}$ is an inadmissible estimator. It is also shown that $\hat{p}$ is generalized Bayes with respect to the uniform prior over $[0, \infty)$. The problems of the confidence intervals for $p$ will be discussed in Section 3. The testing problem for the existence of the signal events will be discussed in Section 4. Section 4 will also discuss a modification of the Type I error rate for future studies on the problem of the existence of signal events.

## 2. (In)admissibility

It is convenient to let $f_\mu$ and $F_\mu$ be the probability mass function (PMF) and the cumulative distribution function (CDF) of the Poisson distribution with mean $\mu$. Then

(1) $$\hat{p}(n) = P(B = n|B \leq n) = \frac{f_b(n)}{F_b(n)}.$$

It is also convenient to abbreviate $\hat{p}(X)$ and $p_{\hat{\theta}}(X)$ by $\hat{p}$ and $p_{\hat{\theta}}$ respectively. Let $\tilde{p}$ be any estimator of $p$. Then under the squared error loss, the risk of $\tilde{p}$ is

$$R(\tilde{p}, \theta) = E_\theta[(\tilde{p} - p)^2] = \sum_{n=0}^{\infty} [\tilde{p}(n) - \frac{b^n}{(b+\theta)^n}]^2 f_{b+\theta}(n). \tag{2}$$

We say $\tilde{p}$ is *inadmissible* if there is an estimator $\tilde{p}'$ for which $R(\tilde{p}', p) \leq R(\tilde{p}, p)$ for all $\theta \geq 0$ with strict inequality for some $\theta$. Otherwise, we say $\tilde{p}$ is *admissible*.

*Inadmissibility of the MLE.* If $\pi$ is a $\sigma$-finite measure over $[0, \infty)$, write $E^\pi$ for integration with respect to the joint distribution of $\theta$ and $X$ when $\theta \sim \pi$; and write $E^\pi(\cdot|n)$ for conditional expectation given $X = n$. Further, let $\bar{R}(\tilde{p}, \pi) = \int_0^\infty R(\tilde{p}, \theta)\pi(d\theta) = E^\pi[(\tilde{p} - p)^2]$ for an estimator $\tilde{p}$ of $p$; and let $\bar{R}(\pi) = \inf_{\tilde{p}} \bar{R}(\tilde{p}, \pi)$. The inadmissibility of the MLE $p_{\hat{\theta}}$ will be obtained from Stein [11].

**Theorem 1.** *If $b \geq 1$, then the MLE $p_{\hat{\theta}} = b^n/(b \vee n)^n$ is inadmissible.*

*Proof.* From Stein [11], a necessary and sufficient condition for $p_{\hat{\theta}}$ to be admissible is that, for every $\theta_0 \in [0, \infty)$ and $\epsilon > 0$, there is a finite prior $\pi$ for which $\pi\{\theta_0\} \geq 1$ and $\bar{R}(p_{\hat{\theta}}, \pi) - \bar{R}(\pi) \leq \epsilon$. In particular, if $p_{\hat{\theta}}$ is to be admissible, then there must be a sequence $\pi_k$ of finite priors for which $\pi_k\{1\} \geq 1$ for very $k = 1, 2, \ldots$ and

$$\lim_{k \to \infty} \bar{R}(p_{\hat{\theta}}, \pi_k) - \bar{R}(\pi_k) = 0. \tag{3}$$

Suppose (3) is true. Let $p_k(n) = E^{\pi_k}[b^n/(b+\theta)^n|n]$. Then

$$\bar{R}(p_{\hat{\theta}}, \pi_k) - \bar{R}(\pi_k) = \sum_{n=0}^{\infty} [p_{\hat{\theta}}(n) - p_k(n)]^2 \int_0^\infty \frac{(b+\theta)^n}{n!} e^{-(b+\theta)} \pi_k(d\theta)$$

$$\geq e^{-(b+1)} \sum_{n=0}^{\infty} \frac{(b+1)^n}{n!} [p_{\hat{\theta}}(n) - p_k(n)]^2.$$

So, $\lim_{k \to \infty} p_k(n) = p_{\hat{\theta}}(n) = b^n/(b \vee n)^n$ for all $n = 0, 1, 2, \ldots$. Note that

$$p_k(n) = \left[\int_0^\infty f_{b+\theta}(n)\pi_k(d\theta)\right]^{-1} \int_0^\infty \frac{b^n}{(b+\theta)^n} f_{b+\theta}(n)\pi_k(d\theta)$$

$$= \left[\int_0^\infty f_{b+\theta}(n)\pi_k(d\theta)\right]^{-1} \int_0^\infty f_b(n)e^{-\theta}\pi_k(d\theta).$$

Let $\pi_k^*(d\theta) = \frac{e^{-\theta}\pi_k(d\theta)}{\int_0^\infty e^{-\theta}\pi_k(d\theta)}$. Then, we have

$$\lim_{k \to \infty} p_k(n) = \lim_{k \to \infty} \frac{\int_0^\infty f_b(n)e^{-\theta}\pi_k(d\theta)}{\int_0^\infty f_{b+\theta}(n)\pi_k(d\theta)} = \lim_{k \to \infty} \frac{\int_0^\infty b^n \pi_k^*(d\theta)}{\int_0^\infty (b+\theta)^n \pi_k^*(d\theta)} = \frac{b^n}{(b \vee n)^n}.$$

Let $\mu_{b,k,n} = \int_0^\infty (b+\theta)^n \pi_k^*(d\theta)$. Note that $\pi_k^*$ is a probability measure. Then, the above expression implies

$$\lim_{k \to \infty} \mu_{b,k,n} = (b \vee n)^n$$

for all $n = 0, 1, 2, \ldots$. This requires

$$\lim_{k \to \infty} \mu_{b,k,n} = \begin{cases} b^n, & \text{if } n \leq m \\ n^n & \text{if } n > m \end{cases} \tag{4}$$

where $m = \lfloor b \rfloor$ is the greatest integer not greater than $b$. Suppose a probability distribution $G$ has moments $\mu_n = b^n$ or $\mu_n = n^n$ with respect to $n \leq m$ or $n > m$ respectively. Then $G$ is determined by its moments since $\sum_{n=0}^{\infty} n^n r^n/n!$ converges if $r < e^{-1}$. Thus (4) requires that the distribution function $G_k(z) = \pi_k^*\{\theta : b+\theta \leq z\}$ converge weakly to the distribution function $G$ (see [5], Section 30). For $b \geq 1$, it requires $\mu_1 = b$. Then $G$ must be degenerate at $b$ so all $\mu_n = b$ for $n \geq 1$. Therefore, such a $G$ can not exist. So, the assumed existence of $\pi_k$ in (3) leads to a contradiction for $b \geq 1$. □

*Admissibility of $\hat{p}$.* The admissibility of $\hat{p}$ will be deduced from a Bayesian approximation. Consider the priors $\pi_\alpha(d\theta) = e^{-\alpha\theta}$ for $\theta \geq 0$ and $\alpha \geq 0$. Then, $\pi_\alpha$ is finite if $\alpha > 0$ and $\pi_0$ is the infinite uniform distribution over $[0, \infty)$. Write $E^\alpha$ for integration with respect to $\pi_\alpha$, $E^\alpha(\cdot|n)$ for posterior expectation, and $\bar{R}(\tilde{p}, \alpha) = E^\alpha[(\tilde{p} - p)^2] = \int_0^\infty R(\tilde{p}, \theta) e^{-\alpha\theta} d\theta$ for the integrated risk of an estimator $\tilde{p}$. This is minimized by $\hat{p}_\alpha = E^\alpha(p|n)$.

The computation of $\hat{p}_\alpha$ is straightforward. First the joint distribution of $(X, \theta)$ under $\pi_\alpha$ is

$$P^\alpha[X = n, d\theta] = \frac{(b+\theta)^n}{n!} e^{-b-(1+\alpha)\theta} d\theta = e^{-b} \sum_{k=0}^{n} \frac{b^k \theta^{n-k}}{k!(n-k)!} e^{-(1+\alpha)\theta} d\theta.$$

The marginal PMF of $X$ under $\pi_\alpha$ is

$$P^\alpha[X = n] = e^{-b} \sum_{k=0}^{n} \frac{b^k}{k!(n-k)!} \int_0^\infty \theta^{n-k} e^{-(1+\alpha)\theta} d\theta$$

$$= e^{-b} \sum_{k=0}^{n} \frac{b^k}{k!(1+\alpha)^{n-k+1}} = \frac{e^{\alpha b}}{(1+\alpha)^{n+1}} F_{(1+\alpha)b}(n).$$

The posterior density of $\theta$ under $\pi_\alpha$ is

$$q_\alpha(\theta) = \left[\frac{e^{\alpha b}}{(1+\alpha)^{n+1}} F_{(1+\alpha)b}(n)\right]^{-1} \frac{(b+\theta)^n}{n!} e^{-b-(1+\alpha)\theta} d\theta.$$

Then

$$\hat{p}_\alpha(n) = E^\alpha(p|n)$$

(5)
$$= \left[\frac{e^{\alpha b}}{(1+\alpha)^{n+1}} F_{(1+\alpha)b}(n)\right]^{-1} \int_0^\infty \frac{b^n}{(b+\theta)^n} \frac{(b+\theta)^n}{n!} e^{-b-(1+\alpha)\theta} d\theta$$

$$= \left[\frac{e^{\alpha b}}{(1+\alpha)^{n+1}} F_{(1+\alpha)b}(n)\right]^{-1} \frac{b^n}{n!(1+\alpha)} e^{-b} = \frac{f_{(1+\alpha)b}(n)}{F_{(1+\alpha)b}(n)}.$$

From the expression of $\hat{p}$ in (1), one has $\hat{p}_0 = \hat{p}$. Let $\delta_{\alpha,n} = f_{(1+\alpha)b}(n)/F_{(1+\alpha)b}(n)$ and denote $\delta_{0,n}$ by $\delta_n$. Then, $\hat{p}_\alpha = \delta_{\alpha,n}$ and $\hat{p} = \delta_n$.

**Theorem 2.** *If $\hat{p}_\alpha$ and $\hat{p}$ are as in (5) and (1), then $\lim_{\alpha \to 0} E^\alpha(\hat{p}_\alpha - \hat{p})^2 = 0$.*

*Proof.* Consider

$$E^\alpha(\hat{p}_\alpha - \hat{p})^2 = \int_0^\infty E_\theta(\hat{p}_\alpha - \hat{p}) e^{-\alpha\theta} d\theta$$

$$= \int_0^\infty \sum_{n=0}^\infty (\delta_{\alpha,n} - \delta_n)^2 \frac{(b+\theta)^n}{n!} e^{-(b+\theta)} e^{-\alpha\theta} d\theta$$

$$\leq \sum_{n=0}^\infty (\delta_{\alpha,n} - \delta_n)^2.$$

Since

$$|\delta_{\alpha,n} - \delta_n|^2 \leq \delta_{\alpha,n}^2 + \delta_n^2 \leq \delta_{\alpha,n} + \delta_n \leq e^{(1+\alpha)b} f_{(1+\alpha)b}(n) + e^b f_b(n)$$

and the right side is summable over $n \geq 0$,

$$\lim_{\alpha \to 0} E^\alpha (\hat{p}_\alpha - \hat{p})^2 \leq \lim_{\alpha \to 0} \sum_{n=0}^\infty (\delta_{\alpha,n} - \delta_n)^2 = \sum_{n=0}^\infty \lim_{\alpha \to 0} (\delta_{\alpha,n} - \delta_n)^2 = 0.$$

from the Dominated Convergence Theorem. $\square$

Since $\hat{p}_\alpha, p$ are less than 1 and for any estimator $\tilde{p}$ with $0 \leq \tilde{p} \leq 1$, $(\hat{p}_\alpha - \tilde{p})(\hat{p}_\alpha - p)$ is absolutely bounded by 4. By the Dominated Convergence Theorem for the second equality below, we have

$$\begin{aligned}
&E^\alpha[(\hat{p}_\alpha - \tilde{p})(\hat{p}_\alpha - p)] \\
&= \int_0^\infty \sum_{n=0}^\infty [\hat{p}_\alpha(n) - \tilde{p}(n)][\hat{p}_\alpha(n) - \frac{b^n}{(b+\theta)^n}] \frac{(b+\theta)^n}{n!} e^{-(b+\theta)} e^{-\alpha\theta} d\theta \\
&= \sum_{n=0}^\infty [\hat{p}_\alpha(n) - \tilde{p}(n)] \int_0^\infty [\hat{p}_\alpha(n) - \frac{b^n}{(b+\theta)^n}] \frac{(b+\theta)^n}{n!} e^{-(b+\theta)} e^{-\alpha\theta} d\theta = 0.
\end{aligned}$$

Then for $\alpha > 0$, $\bar{R}(\tilde{p}, \alpha) = \bar{R}(\hat{p}_\alpha, \alpha) + E^\alpha[(\tilde{p} - \hat{p}_\alpha)^2]$.

**Corollary 1.** $\lim_{\alpha \to 0}[\bar{R}(\hat{p}, \alpha) - \bar{R}(\hat{p}_\alpha, \alpha)] = 0$.

*Proof.* For $\alpha > 0$, $\bar{R}(\hat{p}, \alpha) = E^\alpha[(\hat{p}_\alpha - p)^2 + (\hat{p}_\alpha - \hat{p})^2] = \bar{R}(\hat{p}, \alpha) + E^\alpha[(\hat{p}_\alpha - \hat{p})^2]$ and $E^\alpha[(\hat{p} - \hat{p}_\alpha)^2] \to 0$ as $\alpha \to 0$ by the theorem. $\square$

Since $0 \leq \tilde{p} \leq 1$ for any estimator of $p$ (under consideration), it is clear that from the Dominated Convergence Theorem that $R(\tilde{p}, \theta)$ is continuous in $\theta$ for any estimator $\tilde{p}$ of $p$.

**Corollary 2.** $\hat{p}$ is admissible.

*Proof.* If $\hat{p}$ were inadmissible, then there would be a $\tilde{p}$ for which $R(\tilde{p}, \theta) \leq R(\hat{p}, \theta)$ for all $\theta \geq 0$ and $R(\tilde{p}, \theta_0) < R(\hat{p}, \theta_0)$ for some $\theta_0 \geq 0$. Let $\epsilon_0 = [R(\hat{p}, \theta_0) - R(\tilde{p}, \theta_0)]/2$. Then, there exist an $\eta > 0$, such that $R(\hat{p}, \theta) \geq R(\tilde{p}, \theta) + \epsilon_0$ for all non-negative $\theta$ such that $|\theta - \theta_0| < \eta$. Then

$$\begin{aligned}
\bar{R}(\hat{p}, \alpha) - \bar{R}(\hat{p}_\alpha, \alpha) &\geq \bar{R}(\hat{p}, \alpha) - \bar{R}(\tilde{p}, \alpha) \\
&\geq \int_{\theta_0}^{\theta_0 + \eta} [R(\hat{p}, \theta) - R(\tilde{p}, \theta)] e^{-\alpha\theta} d\theta \\
&\geq \epsilon_0 \frac{e^{-\alpha\theta_0} - e^{-\alpha(\theta_0 + \eta)}}{\alpha} \to \epsilon_0 \eta > 0
\end{aligned}$$

as $\alpha \to 0$, contradicting Corollary 1. $\square$

Of course, it does not follow that the risk of $\hat{p}$ dominates that of $p_{\hat{\theta}}$. The mean squared errors of the two estimators are shown in Figure 1. The graphs show that the MSEs for $\hat{p}$ are greater than the MESs for $p_{\hat{\theta}}$ for small $\theta$ but less than the MSEs for $p_{\hat{\theta}}$ for moderate and large $\theta$.

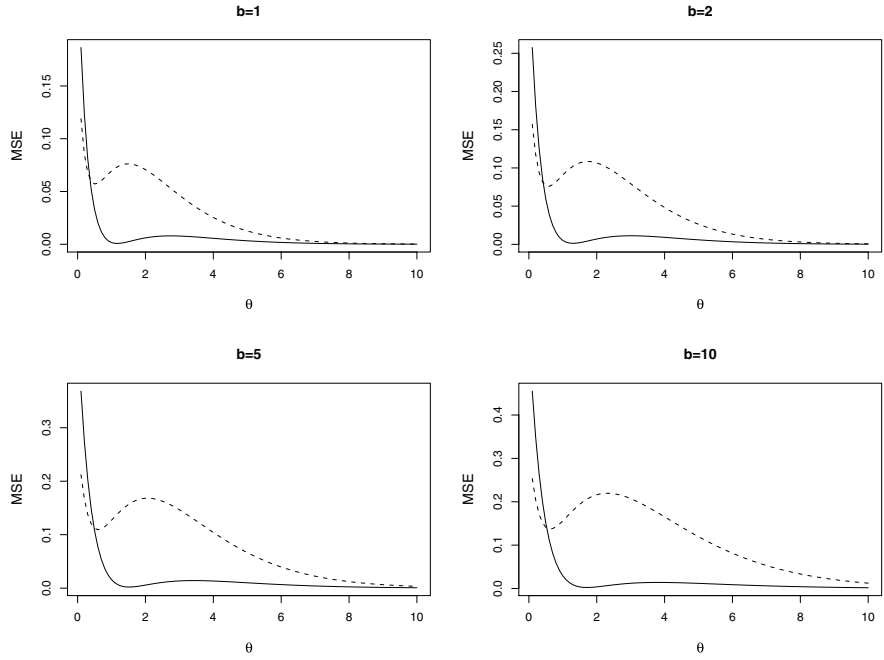

FIG 1. *The MSE for $\hat{p}$ (Solid) and $p_{\hat{\theta}}$ (dotted) for selected b.*

## 3. Confidence and credible interval

Three confidence or credible intervals are considered in this section. They are the unified confidence interval proposed by Feldman and Cousins [8], the conditional frequentist confidence interval proposed by Roe and Woodroofe [10], and the Bayesian credible interval.

The unified confidence interval for $p$ is obtained from the unified approach. It is a transformed interval from the unified confidence interval for $\theta$. In the unified approach, one considers the likelihood ratio

$$R(\theta|k) = \frac{f_{b+\theta}(k)}{f_{b+\hat{\theta}}(k)}.$$

The unified confidence interval for $\theta$ is consisted of taking those $\theta$ for which $R(\theta|k) \geq c(\theta)$, where $c(\theta)$ is the largest value of $c$ for which

$$\sum_{k:R(\theta|k)\geq c} f_{b+\theta}(k) \geq 1 - \alpha.$$

The conditional frequentist confidence interval for $p$ is also a transformed interval from the conditional frequentist interval for $\theta$. In this approach, one modifies the Poisson PMF by the conditional distribution on $B \leq n$ as

$$q_{b,\theta}^n(k) = \begin{cases} f_{b+\theta}(k)/F_b(n), & \text{if } k \leq n \\ \sum_{j=1}^n f_b(j)f_\theta(k-j)/F_b(n), & \text{if } k > n, \end{cases}$$

and then considers the modified likelihood ratio as

$$\tilde{R}^n(\theta|k) = \frac{q_{b,\theta}^n(k)}{\max_{\theta'} q_{b,\theta'}^n(k)}.$$

Let $\tilde{c}_n(\theta)$ be the largest value of $c$ for which

$$\sum_{k:\tilde{R}^n(\theta|k)\geq c} q_{b,\theta}^n(k) \geq 1-\alpha.$$

The conditional frequentist confidence interval consists of those $\theta$ for which $\tilde{R}^n(\theta|k) \geq \tilde{c}_n(\theta)$.

Let $[\tilde{\ell}_u(n), \tilde{u}_u(n)]$ and $[\tilde{\ell}_c(n), \tilde{u}_c(n)]$ be the unified confidence interval and the conditional frequentist interval for $\theta$ respectively. Since both of them are transformation invariant by their definitions, the unified confidence interval $[\ell_u(n), u_u(n)]$ and the conditional frequentist confidence interval $[\ell_c(n), u_c(n)]$ for $p$ are

$$[\ell_u(n), u_u(n)] = [\frac{b^n}{(b+\tilde{u}_u(n))^n}, \frac{b^n}{(b+\tilde{\ell}_u(n))^n}]$$

and

$$[\ell_c(n), u_c(n)] = [\frac{b^n}{(b+\tilde{u}_c(n))^n}, \frac{b^n}{(b+\tilde{\ell}_c(n))^n}]$$

respectively.

Recall that the point estimation $\hat{p}$ of $p$ is also the Bayesian estimator of $p$ under the uniform prior. The Bayesian credible interval for $p$ is considered under the uniform prior for $\theta$ with the least interval length. The method to compute such an interval is introduced in Berger's book ([4], p. 266).

Straightforwardly, the posterior density of $\theta$ under the uniform prior for $\theta$ over $[0, \infty)$ is

$$\tilde{g}(\theta|n) = \frac{f_{b+\theta}(n)}{\int_0^\infty f_{b+\theta}(n)d\theta} = \frac{1}{F_b(n)} \frac{(b+\theta)^n}{n!} e^{-(b+\theta)}.$$

Substituting $\theta = b(p^{-\frac{1}{n}} - 1)$ above and multiplying the absolute value of $d\theta/dp = -bp^{-\frac{n+1}{n}}/n$ on the right side of the above equation, one obtains the posterior density function of $p$, $0 < p \leq 1$, given $X = n$ as

$$g(p|n) = \tilde{g}(b(p^{-1/n} - 1)|n)|\frac{d\theta}{dp}| = \frac{b^{n+1}}{n!nF_b(n)} p^{-\frac{2n+1}{n}} e^{-bp^{-\frac{1}{n}}}.$$

Note that $g(p|n)$ is monotone increasing in $p$ if $b \geq 2n+1$ and is unimodal in $p$ with mode $(\frac{b}{2n+1})^n$ if $b < 2n+1$. Let $p_{\max} = (\frac{b}{2n+1})^n \wedge 1$. Then $g(p|n)$ is uniquely maximized at $p_{\max}$, increasing in $p$ if $p < p_{\max}$ and decreasing in $p$ if $p > p_{\max}$ for $0 < p \leq 1$. The $1-\alpha$ level Bayesian credible interval $[\ell, u] = [\ell(n), u(n)]$ for $p$ can be uniquely solved from

(6) $$G(u|n) - G(l|n) = 1 - \alpha$$

and

(7) $$[\ell, u] = \{p : g(p|n) \geq c_n\}$$

for some $c_n > 0$, where $G(\cdot|n)$ is the posterior distribution of $p$ given $X = n$ under the uniform prior for $\theta$. That the interval solved from (6) and (7) is the shortest interval among those satisfying condition (6) is based on Berger's book.

The algorithm to compute the $1 - \alpha$ level Bayesian credible interval $[\ell, u]$ for $p$ can be described as follows:

TABLE 1
*90% confidence and credible intervals for p when n is close to b for selected b*

| $b$ | $n$ | Unified | Conditional Frequentist | Bayesian |
|---|---|---|---|---|
| 1 | 0 | $[1, 1]$ | $[1, 1]$ | $[1, 1]$ |
| 1 | 1 | $[2.21(10^{-1}), 1]$ | $[2.29(10^{-1}), 1]$ | $[1.62(10^{-1}), 0.862]$ |
| 1 | 2 | $[2.83(10^{-2}), 1]$ | $[2.86(10^{-2}), 1]$ | $[6.78(10^{-3}), 0.495]$ |
| 2 | 1 | $[3.89(10^{-1}), 1]$ | $[4.42(10^{-1}), 1]$ | $[4.00(10^{-1}), 1]$ |
| 2 | 2 | $[1.08(10^{-1}), 1]$ | $[1.15(10^{-1}), 1]$ | $[4.49(10^{-2}), 0.801]$ |
| 2 | 3 | $[2.02(10^{-2}), 1]$ | $[1.95(10^{-2}), 1]$ | $[2.23(10^{-3}), 0.555]$ |
| 2 | 4 | $[2.84(10^{-3}), 0.952]$ | $[2.93(10^{-3}), 1]$ | $[8.21(10^{-5}), 0.273]$ |
| 5 | 4 | $[7.20(10^{-2}), 1]$ | $[1.14(10^{-1}), 1]$ | $[1.53(10^{-2}), 0.823]$ |
| 5 | 5 | $[2.24(10^{-2}), 1]$ | $[3.15(10^{-2}), 1]$ | $[1.06(10^{-3}), 0.712]$ |
| 5 | 6 | $[5.67(10^{-3}), 1]$ | $[6.86(10^{-3}), 1]$ | $[5.85(10^{-5}), 0.564]$ |
| 5 | 7 | $[1.37(10^{-3}), 1]$ | $[1.61(10^{-3}), 1]$ | $[2.78(10^{-6}), 0.384]$ |
| 10 | 9 | $[1.17(10^{-2}), 1]$ | $[2.18(10^{-2}), 1]$ | $[3.33(10^{-5}), 0.719]$ |
| 10 | 10 | $[4.02(10^{-3}), 1]$ | $[6.68(10^{-3}), 1]$ | $[2.10(10^{-6}), 0.633]$ |
| 10 | 11 | $[1.53(10^{-3}), 1]$ | $[1.75(10^{-3}), 1]$ | $[1.20(10^{-7}), 0.528]$ |
| 10 | 12 | $[4.07(10^{-4}), 1]$ | $[4.51(10^{-4}), 1]$ | $[6.29(10^{-9}), 0.406]$ |
| 10 | 13 | $[9.71(10^{-5}), 0.834]$ | $[1.18(10^{-4}), 1]$ | $[3.15(10^{-10}), 0.278]$ |
| 10 | 14 | $[2.09(10^{-5}), 0.501]$ | $[2.22(10^{-5}), 1]$ | $[1.54(10^{-11}), 0.163]$ |

(i) Find a $z$ so that $G(z|n) = \alpha$.
(ii) If $g(z|n) \leq g(1|n)$, then $\ell = z$ and $u = 1$. Otherwise, repeat step (iii) until convergence.
(iii) Let $c_1 = g(z|n)$ and $c_2 = g(1|n)$. Let $c_0 = (c_1 + c_2)/2$. Find $p_\ell < p_{\max} < p_u$ so that $g(p_\ell) = g(p_u) = c_0$. If $G(p_u) - G(p_\ell) > 1 - \alpha$ then let $c_2 = c_0$; otherwise let $c_1 = c_0$ and iterate.

Table 1 lists the lower bound and upper bound of the confidence and credible intervals for $p$ as functions of $n$ when $n$ is close to the background parameter $b$ for selected $b$.

The frequentist coverage probability of the confidence and credible intervals for $p$ is the probability of the interval to contain the true value of $p_\theta(n) = b^n/(b+\theta)^n$. The numerical results of the frequentist coverage probability of the 90% confidence and credible intervals are given in Figure 2. As suggested by Figure 2, the Bayesian credible interval has a low coverage when $\theta$ is close to 0 and it steadily increases to almost 0.9 as $\theta$ becomes large. However, since the length of the Bayesian credible interval is much shorter than the conditional frequentist and the unified intervals from the numerical result, Figure 2 does not suggest that the Bayesian credible interval is worse than the other two. The curves of coverage probabilities fluctuate up and down since the coverage probabilities are discontinuous in $\theta$.

## 4. A testing problem

An interesting problem in particle physics is to determine if there exist any signal events in an experiment. In statistics, this problem can be described as a testing problem for the null hypothesis of $H_0 : S = 0$ versus the alternative hypothesis of $H_1 : S > 0$. Recall that $p = P_\theta[S = 0|X = n] = b^n/(b+\theta)^n$ and note that $\theta$ is unknown. The probability $p$ should be estimated. As we have seen, two estimators have been discussed. The first one is $\hat{p}$ and the second one is the MLE $p_{\hat{\theta}}$. Since $p_{\hat{\theta}}$ is inadmissible, we only discuss $\hat{p}(n) = f_b(n)/F_b(n)$ here.

An estimation problem here is to estimate the indicator function $I_{S>0}$ instead of its probability. As an estimator of $I_{S>0}$, we will show that $1 - \hat{p}$ is also admissible under the squared error loss, where the squared error loss is defined by $L(\delta, S) = (\delta - I_{S>0})^2$ for an estimator $\delta$ of $I_{S>0}$.

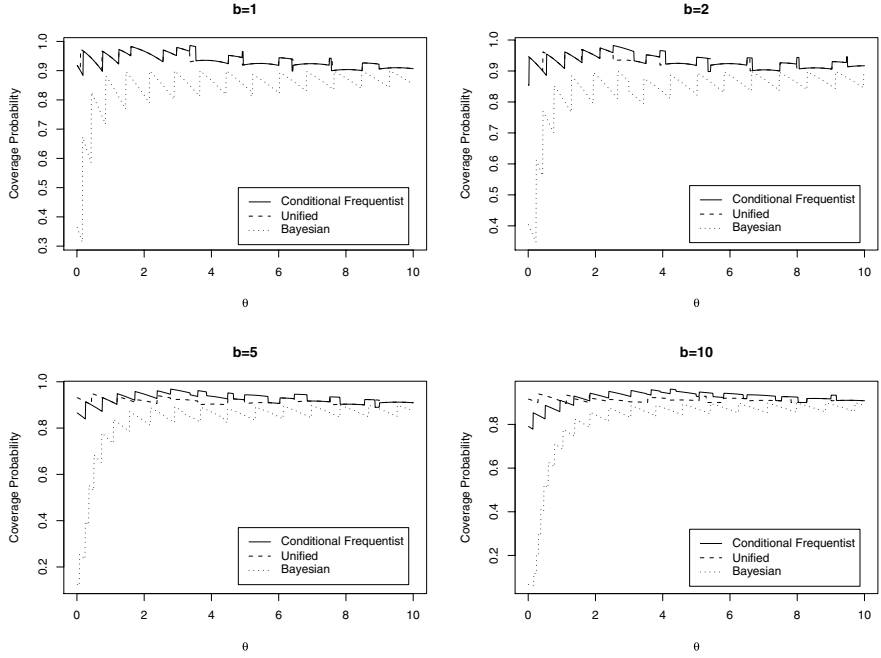

FIG 2. *Frequentist coverage probability of the 90% confidence and credible intervals for p as functions of $\theta$ for selected b.*

As before, we still consider the priors $\pi_\alpha(\theta) = e^{-\alpha\theta}$, $\theta \geq 0$ for $\alpha \geq 0$, but here we define $R(\delta, \theta) = E_\theta[(\delta - I_{S>0})^2]$ and $\bar{R}(\delta, \alpha) = E^\alpha[(\delta - I_{S>0})^2] = \int_0^\infty R(\delta, \theta)e^{-\alpha\theta}d\theta$. This is minimized by $E^\alpha(I_{S>0}|n)$.

The computation of $E^\alpha(I_{S>0}|n)$ is still straightforward. SInce $P_\theta[B = k, X = n] = b^k\theta^{n-k}e^{-(b+\theta)}/[k!(n-k)!]$ follows the joint posterior mass function, for $n \geq 0$,

$$P^\alpha[B = k, X = n] = \int_0^\infty P_\theta[B = k, X = n]e^{-\alpha\theta}d\theta = \frac{b^k e^{-b}}{k!(1+\alpha)^{n-k+1}},$$

and the marginal posterior mass function of X, for $n \geq 0$,

$$P^\alpha[X = n] = \sum_{k=0}^n P^\alpha[B = k, X = n] = \frac{e^{\alpha b}}{(1+\alpha)^{n+1}} F_{(1+\alpha)b}(n).$$

Then, the conditional posterior mass function of $B = k$ given $X = n$ is

$$P^\alpha[B = k|X = n] = \frac{(1+\alpha)^k e^{-(1+\alpha)b}}{k! F_{(1+\alpha)b}(n)} = \frac{f_{(1+\alpha)b}(k)}{F_{(1+\alpha)b}(n)}.$$

This gives

$$E^\alpha(I_{S>0}|n) = \sum_{k=0}^{n-1} P^\alpha[B = k|X = n] = \frac{F_{(1+\alpha)b}(n-1)}{F_{(1+\alpha)b}(n)}.$$

Since $E^\alpha(I_{S>0}|n) = 1 - \hat{p}_\alpha(n)$ where $\hat{p}_\alpha(n)$ is given by (5) and since for $\alpha > 0$ there is also the identity $\bar{R}(\hat{p}, \alpha) = \bar{R}(\hat{p}, \alpha) + E^\alpha[(\hat{p}_\alpha - \hat{p})^2]$, by Theorem 2 one can see that $1-\hat{p}$ is also an admissible estimator of the indicator function $I_{S>0}$. Details of the proof is omitted here since there are many overlaps in the proof.

The mean squared errors of the two estimators $\hat{p}$ and $p_{\hat{\theta}}$ as estimators of $I_{S>0}$ are identical to the mean squared errors of the two estimators as estimators of $1-p$, which can also be seen in Figure 1. However as estimators of $I_{S>0}$, one can also study the conditional risk, where the conditional risk of an estimator $\delta$ of $I_{S>0}$ is defined by $E_\theta[(\delta - I_{S>0})^2 | B \leq n]$. As a result, the conditional risk of $1-\hat{p}$ is

$$(8) \qquad E_\theta[(\hat{p} - I_{S=0})^2 | B \leq n] = V(I_{B=n} | B \leq n) = \frac{f_b(n) F_b(n-1)}{F_b^2(n)}$$

and the conditional risk of the MLE $1 - p_{\hat{\theta}}(n)$ is

$$(9) \qquad \begin{aligned} E_\theta[(p_{\hat{\theta}} - I_{S=0})^2 | B \leq n] &= V(I_{B=n} | B \leq n) + [p_{\hat{\theta}}(n) - \hat{p}(n)]^2 \\ &= \frac{f_b(n) F_b(n-1)}{F_b^2(n)} + \left[\frac{f_b(n)}{F_b(n)} - \frac{b^n}{(b \vee n)^n}\right]^2. \end{aligned}$$

Obviously they do not depend on $\theta$ and (8) is less than (9). The conditional risks of $1 - \hat{p}$ and the MLE $1 - p_{\hat{p}}$ as functions of $n$ for selected $b$ are displayed in Figure 3.

Type I error rate can be modified in the testing problem of the null hypothesis of $H_0 : S = 0$ versus the alternative hypothesis of $H_1 : S > 0$. Suppose the null hypothesis is rejected if $n \geq n_0$, where $n_0$ is a fixed integer. Then, the classical type I error rate is

$$\alpha(n_0) = P_\theta[X \geq n_0 | S = 0] = P_b[X \geq n_0] = 1 - F_b(n_0 - 1).$$

Based on the estimator $\hat{p}$ of $p$, a modified type I error rate is proposed as follows.

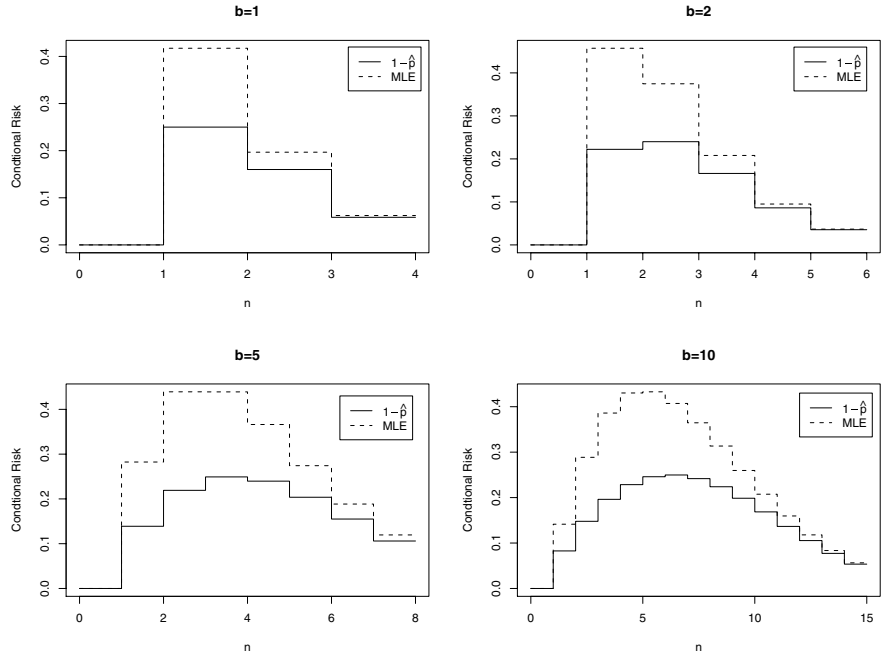

FIG 3. *Conditional risk of estimators of $I_{S>0}$ as functions of $n$ for selected $b$.*

TABLE 2
*Least $n_0$ for $\alpha(n_0) \leq 0.1$ and $\alpha^*(n_0) \leq 0.1$ respectively*

| $b$ | Classical | Modified |
|---|---|---|
| 1 | 3 | 3 |
| 2 | 5 | 5 |
| 3 | 6 | 6 |
| 4 | 8 | 7 |
| 5 | 9 | 8 |
| 6 | 10 | 9 |
| 7 | 11 | 10 |
| 8 | 13 | 11 |
| 9 | 14 | 12 |
| 10 | 15 | 13 |

Note that

$$P_\theta[X \geq n_0 | S = 0] = \frac{P_\theta[X \geq n_0, S = 0]}{P_\theta[S = 0]} = \frac{\sum_{n=n_0}^\infty P_\theta[X = n, S = 0]}{\sum_{n=0}^\infty P_\theta[X = n, S = 0]}$$
$$= \frac{\sum_{n=n_0}^\infty P_\theta[S = 0 | X = n] P_\theta[X = n]}{\sum_{n=0}^\infty P_\theta[S = 0 | X = n] P_\theta[X = n]} = \frac{\sum_{n=n_0}^\infty p_\theta(n) f_{b+\theta}(n)}{\sum_{n=0}^\infty p_\theta(n) f_{b+\theta}(n)}.$$

Substituting $p_\theta(n)$ by $\hat{p}(n)$ in both the numerator and denominator in the last equation, we have

$$\alpha_\theta^*(n_0) = \frac{\sum_{n=n_0}^\infty \hat{p}(n) f_{b+\theta}(n)}{\sum_{n=0}^\infty \hat{p}(n) f_{b+\theta}(n)}.$$

Since $\alpha_\theta^*(n_0)$ is strictly decreasing in $\theta$ for any $n_0 > 0$, a modified type I error rate is proposed by

$$\alpha^*(n_0) = \alpha_0^*(n_0) = \frac{\sum_{n=n_0}^\infty \hat{p}(n) f_b(n)}{\sum_{n=0}^\infty \hat{p}(n) f_b(n)}.$$

Table 2 lists the least $n_0$ for $\alpha(n_0) \leq 0.1$ and $\alpha^*(n_0) \leq 0.1$ respectively for selected $b$. It can be seen that the values of $n_0$ based on modified type I error rates are always no greater than those based on the classical type I error rates.

## 5. Remarks

An intrinsic testing problem in particle physics is to know the general conclusion of the possibility of the signal events in experiments, which can be described as a testing problem for the null hypothesis of $H_0 : \theta = 0$ versus the alternative hypothesis of $H_1 : \theta > 0$. Since $P_\theta[X = n | S = 0] = P_0[X = n]$, the two testing problems are equivalent if the inference is only based on observations of a single experiment.

Conditioning on $B \leq X$ when $X = n$ in the signal plus background model was first proposed by Roe and Woodroofe [10]. Later Woodroofe and Wang [13] consider a testing problem for $H_0 : \theta \geq \theta_0$ versus $H_1 : \theta < \theta_0$ for a positive $\theta_0$ based on the condition $B \leq n$. Later on, Zhang and Woodroofe [14] consider an estimation problem for the signal event still based on the condition $B \leq n$. This paper considers the estimation of the probability of the existence of signal events based on the condition $B \leq n$.

**References**

[1] AGUILAR, A., ET AL., (2001). Evidence for neutrino oscillations from the observation of $\bar{\nu}_e$ appearance in a $\bar{\nu}_\mu$ beam. *Physical Review D* **64** 112007.
[2] ARMBRUSTER, B., ET AL. (2002). Upper limits for neutrino oscillations $\bar{\nu}_\mu \to \bar{\nu}_e$ from muon decay at rest. *Physical Review D* **65** 112001.
[3] AVVAKUMOV, S., ET AL. (2002). Search for $\nu_\mu \to \nu_e$ and $\bar{\nu}_e \to \bar{\nu}_e$ oscillations at NuTev. *Physical Review Letters* **89** 011804.
[4] BERGER, J. (1980). *Statistical Decision Theory, Foundations, Concepts, and Methods*. Spring-Verlag. MR0580664
[5] BILLINGSLEY, P. (1995). *Probability and Measure*. Wiley, New York. MR1324786
[6] EIKEL, K. (2003). The LSND and KARMEN short baseline accelerator-based neutrino oscillation searchs. Conference on the Intersections of Particle and Nuclear Physics, CIPANP2003, May 29–24, 2003, New York City. Available on http://www-ik1.fzk.de/www/karmen.
[7] EITEL, K., ET AL. (2000). Update of the KARMEN2 $\nu_\mu \to \nu_e$ oscillation search. In *Proceedings 14th Lake Louise Winter Institute Electroweak Physics*. World Scientific, Singapore, pp. 353–360.
[8] FELDMAN, G. J. AND COUSINS, R. (1998). Unified approach to the classical statistical analysis of small signals. *Physical Review D* **57** 3873–3889.
[9] MANDELKERN, M. (2002). Setting confidence intervals for bounded parameters (with discussions). *Statistical Science* **17** 136–172. MR1939335
[10] ROE, B. AND WOODROOFE, M. (1999). Improved probability method for estimating signal in the presence of background. *Physical Reviews D* **60** 053009.
[11] STEIN, C. (1955). A necessary and sufficient condition for admissibility. *Ann. Math. Statist.* **26** 518–522. MR0070929
[12] WOODROOFE, M. AND ROE, B. (2003). Statistics issues for the MiniBoone experiment. Technical Report, University of Michigan. Available on http://www.ippp.dur.ac.uk/Workshops/02/statistics/proceedings/roe.ps.
[13] WOODROOFE, M. AND WANG, H. (2000). The problem of low counts in a signal plus noise model. *Annals of Statistics* **28** 1561–1569. MR1835031
[14] ZHANG, T. AND WOODROOFE, M. (2005). Admissible minimax estimation of the signal with known background. *Statistica Sinica* **15** 59–72. MR2125720